\documentclass[12pt]{amsart}
\usepackage{amssymb,amsthm,amsfonts,amsmath}  \newcommand{\Q}{{\mathbb Q}}
 
\newcommand{\Z}{{\mathbb Z}} \newcommand{\N}{{\mathbb N}}
\newcommand{\R}{{\mathbb R}} 
\newcommand{\Proj}{{\mathbf P}} 
 
 \newcommand{\EE}{{\mathcal E}}
 \newcommand{\C}{{\mathbb C}}

\newcommand{\F}{{\mathbb F}} 
\newcommand{\Fpbar}{{\overline{\mathbb{F}}_p}}

 \newcommand{\W}{{\mathcal{W}}}

\newcommand{\ord}{{\operatorname{ord}}}
\newcommand{\Com}{{\operatorname{Com}}}

\newcommand{\comment}[1]{} \newtheorem{theorem}{Theorem}[section]

\newtheorem{prop}[theorem]{Proposition}
 \theoremstyle{definition}
 
\newtheorem{definition}{Definition} \theoremstyle{remark}
\newtheorem{rem}{Remark} \newtheorem{example}{Example} \author{Kirsten
  Eisentr\"ager} \title[Hilbert's Tenth Problem for function
fields]{Hilbert's Tenth Problem for function fields of characteristic
  zero}
\begin{document}
\begin{abstract}
  In this article we outline the methods that are used to prove
  undecidability of Hilbert's Tenth Problem for function fields of
  characteristic zero. Following Denef we show how rank one elliptic
  curves can be used to prove undecidability for rational function
  fields over formally real fields. We also sketch the
  undecidability proofs for function fields of varieties over the
  complex numbers of dimension at least 2.
\end{abstract}
\maketitle
\section{Introduction}
Hilbert's Tenth Problem in its original form was to find an algorithm
to decide, given a polynomial equation $f(x_1,\dots,x_n)=0$ with
coefficients in the ring $\Z$ of integers, whether it has a solution
with $x_1,\dots,x_n \in \Z$.  Matijasevi\v{c}~(\cite{Mat70}), based on
work by Davis, Putnam and Robinson~(\cite{DPR61}), proved that no such
algorithm exists, {\it i.e.}\ Hilbert's Tenth Problem is undecidable.
Since then, analogues of this problem have been studied by asking the
same question for polynomial equations with coefficients and solutions
in other commutative rings $R$.  We will refer to this as {\sl
  Hilbert's Tenth Problem over $R$}.  Perhaps the most important
unsolved question in this area is the case $R=\Q$.  There has been
recent progress by Poonen~(\cite{Poo03}) who proved undecidability for
large subrings of $\Q$. The function field analogue, namely Hilbert's
Tenth Problem for the function field $k$ of a curve over a finite
field, is undecidable.  This was proved by Pheidas~(\cite{Phei91}) for
$k=\F_q(t)$with $q$ odd, and by Videla~(\cite{Vi94}) for $\F_q(t)$
with $q$ even.  Shlapentokh~(\cite{Sh2000}) generalized Pheidas'
result to finite extensions of $\F_q(t)$ with $q$ odd and to certain
function fields over possibly infinite constant fields of odd
characteristic, and the remaining cases in characteristic $2$ are
treated in \cite{Eis2002}.  Hilbert's Tenth Problem is also known to
be undecidable for several rational function fields of characteristic
zero: In 1978 Denef proved the undecidability of Hilbert's Tenth
Problem for rational function fields $K(t)$ over formally real fields
$K$~(\cite{Den78}). Kim and Roush~(\cite{KR92}) showed that the
problem is undecidable for the purely transcendental function field
$\C(t_1,t_2)$. In \cite{Eis04} this was
generalized to finite extensions of $\C(t_1,\dots,t_n)$ with $n \geq 2$.
The problem is also known to be undecidable for function fields over
$p$-adic fields (\cite{KR95}, \cite{Eispadic}, \cite{MB}), and for
function fields over formally real fields (\cite{MB}).  In
Hilbert's Tenth Problem the coefficients of the equations have to be
input into a Turing machine, so when we consider the problem for
uncountable rings we restrict the coefficients to a subring $R'$ of
$R$ which is finitely generated as a $\Z$-algebra.  We say that {\em
  Hilbert's Tenth Problem for $R$ with coefficients in $R'$} is
undecidable if there is no algorithm that decides whether or not
multivariate polynomial equations with coefficients in $R'$ have a
solution in $R$.  In this paper we will discuss the undecidability
proofs for $\R(t)$, and $\C(t_1,t_2)$. In these cases, we consider
polynomials with coefficients in $\Z[t]$ and $\Z[t_1,t_2]$,
respectively.

The biggest open problems for function fields are Hilbert's Tenth
Problem for $\C(t)$ and for $\Fpbar(t)$.

\section{Approach for function fields of characteristic zero}
\subsection{Preliminaries}
Before we can describe the approach that is used in characteristic
zero we need the following definition. All the rings we consider are
commutative with $1$.
\begin{definition}
Let $R$ be a ring. A subset $Q \subseteq R^k$ is {\em diophantine
over $R$} if there exists a polynomial $f(x_1, \dots,
x_k,y_1,\dots,y_m) \in R[x_1,\dots,x_k,$ $y_1,\dots,y_m]$ such that
\[
Q = \{ \vec x \in R^k: \exists \, y_1, \dots, y_m \in R: f(\vec x, y_1,
\dots, y_m) =0\}.
\] 

Let $R'$ be a subring of $R$ and suppose that $f$ can be chosen such
that its coefficients are in $R'$. Then we say that $Q$ is {\em diophantine
over $R$ with coefficients in $R'$}.
\end{definition}

\begin{example}
The set of natural numbers $\N$ is diophantine over $\Z$. This follows
from the fact that every natural number can be written as a sum of
four squares, so 
\[
\N = \{ a \in \Z : \exists \, y_1, \dots, y_4 \in \Z : (y_1 ^2 + y_2^2 + y_3^2 +
y_4^2 -a =0) \}.
\]
\end{example}
\begin{example}
  The set of primes is diophantine over $\Z$. This follows from the
  proof of Hilbert's Tenth Problem for $\Z$, where it was shown that
  every recursively enumerable subset of $\Z$ is diophantine over $\Z$
  (\cite{DPR61,Mat70}). Clearly the prime numbers form a recursively
  enumerable subset of $\Z$.
\end{example}

\subsection{Combining diophantine equations}
If $K$ is a formally real field and $P_1=0$, $P_2=0$ are two diophantine
equations over $K(t)$ with coefficients in $\Z[t]$, then $P_1=0 \land P_2=0$
and $P_1 =0\lor P_2 =0$ are also diophantine with coefficients in $\Z[t]$:
we have
\[ P_1 = 0 \land P_2 =0 \leftrightarrow P_1^2 + tP_2^2 =0
\]
and 
\[
P_1= 0 \lor P_2 =0 \leftrightarrow P_1 P_2 =0.
\]
The same holds for diophantine equations over $\C(t_1, t_2)$ with
coefficients in $\Z[t_1,t_2]$:
\[
P_1 =0 \land P_2=0 \leftrightarrow P_1^2 + t_1P_2^2 =0 \text{ and }P_1 =0 \lor P_2 =0 \leftrightarrow P_1P_2=0.
\]

An argument similar to the one above can be made for many other rings
whose quotient field is not algebraically closed, provided that the
ring of coefficients is large enough.

The above argument shows that proving undecidability of Hilbert's
Tenth Problem for $K(t)$ with coefficients in $\Z[t]$ is the same as
proving that the positive existential theory of $K(t)$ in the language
$\langle \, +,\cdot\;;0,1,t \rangle$ is undecidable. Similarly, Hilbert's Tenth
Problem for $\C(t_1,t_2)$ with coefficients in $\Z[t_1,t_2]$ is
undecidable if and only if the positive existential theory of
$\C(t_1,t_2)$ in the language $\langle\, +,\cdot\,;0,1,t_1,t_2\rangle$ is
undecidable.
\subsection{Approach in characteristic zero}
We can use a reduction argument to prove undecidability for a ring $R$
of characteristic zero if we can give a diophantine definition of $\Z$
inside $R$.  We have the following proposition:
\begin{prop}\label{diophantinedefinition}
  Let $R$ be an integral domain of characteristic zero. Let $R'$ be a
  subring of $R$, which is finitely generated as a $\Z$-algebra and
  such that the fraction field of $R$ does not contain an algebraic
  closure of $R'$. Assume that $\Z$ is a diophantine subset of $R$
  with coefficients in $R'$. Then Hilbert's Tenth Problem for $R$ with
  coefficients in $R'$ is undecidable.
\end{prop}
\begin{proof}
  Given a polynomial equation $f(x_1,\dots,x_n)=0$ over $\Z$ we can
  construct a system of polynomial equations over $R$ with
  coefficients in $R'$ by taking the original equation together with,
  for each $i=1,\dots,n$ an equation $g_i(x_i,\dots)=0$ involving
  $x_i$ and a new set of variables, such that in any solution over $R$
  of the system, $g_i=0$ forces $x_i$ to be in $\Z$. In other words,
  the new system of equations has a solution over $R$ if and only
  $f(x_1,\dots,x_n)$ has a solution in $\Z$. Also, since the quotient
  field of $R$ does not contain the algebraic closure of $R'$, the
  system over $R$ with coefficients in $R'$ is equivalent to a single
  polynomial equation with coefficients in $R'$ (\cite[p.\ 51]{PheZa00}).
\end{proof}

Sometimes we cannot give a diophantine definition of the integers
inside a ring $R$, but we might be able to construct a model of the
integers inside $R$.

\begin{definition} A {\em diophantine model} of $\langle\Z,0,1;+,\cdot \rangle$ over
  $R$ is a diophantine subset $S \subseteq R^m$ equipped with a bijection
  $\phi: \Z \to S$ such that under $\phi$, the graphs of addition and
  multiplication correspond to diophantine subsets of $S^3$.

  Let $R'$ be a subring of $R$. A {\em diophantine model} of
  $\langle\Z,0,1;+,\cdot \rangle$ over $R$ {\em with coefficients in $R'$} is a
  diophantine model of $\langle\Z,0,1;+,\cdot \rangle$, where in addition $S$ and
  the graphs of addition and multiplication are diophantine over $R$
  with coefficients in $R'$.
\end{definition}
A similar argument as for Proposition~\ref{diophantinedefinition} can
be used to prove the following
\begin{prop}
  Let $R,R'$ be as in Proposition~\ref{diophantinedefinition}. If we
  have a diophantine model of $\langle \Z,0,1;+ , \cdot \rangle $ over $R$ with
  coefficients in $R'$, then
  Hilbert's Tenth Problem for $R$ with coefficients in $R'$ is undecidable.
\end{prop}

\section{Fields of rational functions over formally real fields}

In this section we will give an outline of Denef's theorem:
\begin{theorem}[\cite{Den78}]\label{thm1}
Let $K$ be a formally real field, i.e.\ $-1$ is not the sum of
squares. Hilbert's Tenth Problem for $K(t)$ with
coefficients in $\Z[t]$ is undecidable.
\end{theorem}
We follow Denef's proof and use elliptic curves to construct a model
of $\langle \Z,0,1;+ , \cdot \rangle $ in $K(t)$. Denef actually gives a
diophantine definition of $\Z$ inside $K(t)$, but we will construct a
model of the integers, because this approach is slightly shorter and it is
used in subsequent papers (\cite{KR92},\cite{KR95}, \cite{Eispadic}).

To construct the diophantine model of the integers we first have to
obtain a set $S$ which is diophantine over $K(t)$ and which has a natural
bijection to $\Z$.

\subsection{Obtaining $S$}\label{obtainS}
Let $E$ be an elliptic curve defined over $\Q$ without complex
multiplication and with Weierstrass equation
\begin{equation}\label{curve}
  y^2 =x^3 + ax + b.\end{equation} The points $(x,y)$ satisfying
equation~(\ref{curve}) together with the
point at infinity form an abelian group, and the group law is given by
equations.  We will now look at a twist of $E$, the elliptic curve
$\EE$, which is defined to be the smooth projective model of
\[
(t^3+at+b) Y^2 = X^3 + aX+b. \label{eq1}
\]
This is an elliptic curve defined over the rational function field $\Q(t)$. An
obvious point on $\EE$ which is defined over $\Q(t)$ is the point
$P_1:=(t,1)$. Denef proved:
\begin{theorem}\cite[p.\ 396]{Den78}\label{pointorder}
  The point $P_1$ has infinite order and generates the group
  $\EE(K(t))$ modulo points of order 2.
\end{theorem}

The elliptic curve $\EE$ is a projective variety, but any
  projective algebraic set can be partitioned into finitely many
  affine algebraic sets, which can then be embedded into a single
  affine algebraic set.  This implies that the set
  \[\EE(K(t))=\{(X,Y):X,Y \in K(t) \land (t^3 + at+b)Y^2 = X^3 + aX+b\} \cup
  \{ \mathbf{O}\}\]is diophantine over $K(t)$, since we can take care of the
  point at infinity $\mathbf{O}$ of $\EE$.

  We can express by polynomial equations that a point on the elliptic
  curve is of the form $2 \cdot P$, where $P$ is another point on the
  curve.  Hence the set
  \begin{align*}{S'}=&\{(X_{2n},Y_{2n}):n \in \Z\}\\=&\{(x,y)
\in (K(t))^2: \\&\exists \,u,v \in K(t):(u,v) \in \EE(K(t)) \land (x,y)=2(u,v)\}
\end{align*} is  diophantine over $K(t)$ with coefficients in $\Z[t]$. 
Then the set 
 \begin{align*}{S''}&=\{(X_{n},Y_{n}):n \in \Z\}\\&=\{(x,y)
   \in (K(t))^2: \exists \,n \in \Z:\\ &\left( (x,y)=(X_{2n},Y_{2n})
     \lor(x,y)=(X_{2n},Y_{2n})+P_1\right)\}
\end{align*} is diophantine over $K(t)$ with coefficients in $\Z[t]$ as well.

Let $P_n:=n \cdot (t,1) = (X_n,Y_n)$ for $n \in \Z - \{0\}$, and let
$P_0:= \mathbf{O}$.  The set $S'
\cup S''$ is equal to the set $\{P_n: n \in \Z\}$.
Let $Z_n:=\frac{X_n}{tY_n}$ for $n \in \Z - \{0\}$, and let $Z_0:=0$.
We define $S$ to be the set \[S = \left\{ Z_{n}: n \in \Z \right\}.\]

Then $S$ is diophantine over $K(t)$.
Since $Z_n \in \Q(t)$, we can consider $Z_n$ as a function on
the projective line $\Proj^1_{\Q}= \Q \cup \{\infty\}$. Denef (\cite[p.\
396]{Den78}) proved the following proposition:
\begin{prop}\label{denef}
  Considered as a function on $\Proj^1_{\Q}$, $Z_n$ takes the value
  $n$ at infinity.
\end{prop}

For $n \neq m$, we have $Z_n \neq Z_m$, and so by associating
the point $Z_n$ to an integer $n$ we obtain an obvious bijection
between $\Z$ and ${S}$. This is the set that we will use for the
diophantine model of $\langle\Z,0,1;+,\cdot \rangle$.
\subsection{Existentially defining multiplication and addition}
The bijection $\phi: \Z \to S$ given by $\phi(n) = Z_n$ induces
multiplication and addition laws on the set $S$, and it remains to
show that the graphs of addition and multiplication on $S$ are
diophantine over $K(t)$.  This means that we have to show that the
sets
\[
S_{add}:=\{(Z_n,Z_m,Z_{\ell}) \in S^3: n+m =\ell \}
\]
and
\[
S_{mult}:= \{(Z_n,Z_m,Z_{\ell}) \in S^3: n \cdot m =\ell \}
\]
are diophantine over $K(t)$.  Since addition of points on the elliptic
curve is given by equations involving rational functions of the
coordinates of the points, it follows easily that the set $S_{add}$ is
diophantine over $K(t)$:
\begin{gather*}
Z_n + Z_m = Z_{\ell} \leftrightarrow \exists (X,Y), (X',Y'), (X'',Y'') \in
\EE(K(t)):\\ \left(Z_n = \frac{X}{tY}, Z_m=\frac{X'}{tY'},
Z_{\ell}=\frac{X''}{tY''} \land (X,Y)+(X',Y') =(X'',Y'')\right).
\end{gather*}

The difficult part is showing that $S_{mult}$ is diophantine.
\subsection{Defining multiplication}
We define the discrete valuation \linebreak $\ord_{t^{-1}}: K(t) \to
\Z\, \cup\, \{\infty\}$ by $\ord_{t^{-1}}u = - \deg f + \deg g$, for $u \in
K(t)^{*}$, $u = f/g$ with $f,g \in K[t]$. We let $\ord_{t^{-1}}(0)=\infty$.
Proposition~\ref{denef} implies that for $n \neq 0$,
$\ord_{t^{-1}}(Z_n)=0$ and $\ord_{t^{-1}}(Z_n-n)>0$.

The discrete valuation $\ord_{t^{-1}}$ has the following properties:
For $u\in K(t)$ we have
\begin{enumerate}
\item $\ord_{t^{-1}}(u)=0$ if and only if $u$ takes a
nonzero value $a \in K$ at infinity.
\item $\ord_{t^{-1}}(u) >0$ if and only if $u$ takes the value zero at
  infinity.
\item $\ord_{t^{-1}}(u) < 0$ if and only if $u$ takes the value
  infinity at infinity. 
\end{enumerate}

We will use the discrete valuation $\ord_{t^{-1}}$ to existentially
define multiplication of elements of $S$.
This is done in the following theorem.
\begin{theorem}\label{DenefThm}
  Assume that the set $T':= \{ u \in K(t): \ord_{t^{-1}}(u) >0 \}$ is
  diophantine over $K(t)$. Then the set $S_{mult}$ is diophantine over
  $K(t)$, {\it i.e.\ }we can existentially define multiplication of
  elements of $S$.
\end{theorem}
\begin{proof} The theorem follows immediately from the following
  claim.\\
  {\bf Claim:} Given $n,m, \ell \in \Z$ we have $n \cdot m = \ell$ if and
  only if $Z_{n} \cdot Z_{m} - Z_{\ell} \in T'$, {\it i.e.\ } $$\ord_{t^{-1}}(Z_{n} \cdot
  Z_{m} - Z_{\ell})>0.$$

  {\bf Proof of Claim:} If $n \cdot m = \ell$, then by 
  Proposition~\ref{denef}, $Z_{n} \cdot Z_{m} - Z_{\ell}$ takes the
  value $n \cdot m - \ell =0$ at infinity, and hence
  $\ord_{t^{-1}}(Z_{n} \cdot Z_{m} - Z_{\ell})>0$.

  If $n \cdot m \neq \ell$, then $Z_{n} \cdot Z_{m}-Z_{\ell}$ takes the
  value $n\cdot m-  \ell \neq 0$ at infinity. Hence $\ord_{t^{-1}}(Z_{n} \cdot
  Z_{m} - Z_{\ell})=0$.
This proves the claim.

Since we assumed that the set $T'$ of all elements with positive
valuation at $t^{-1}$ was diophantine over $K(t)$ this proves the
theorem.
\end{proof}
\begin{rem}\label{rem1}
  We can modify the set $T'$ and still make the proof of
  Theorem~\ref{DenefThm} work. What we needed in the proof was a
  diophantine set $T$ with the following properties:
\begin{enumerate}
\item If $Z \in \Q(t)$ and $\ord_{t^{-1}}( Z )>0$, then $Z \in T$.
\item If $Z \in K(t)$
and $Z \in T$, then $\ord_{t^{-1}} (Z) >0$.
\end{enumerate}
This is enough because the functions $Z_n$ are elements of $\Q(t)$.
\end{rem}

\subsection{How to obtain a diophantine definition for $T$}
We will now define the set $T$ that has the properties in
Remark~\ref{rem1}. By Theorem~\ref{DenefThm} this is enough to finish
the proof of Theorem~\ref{thm1}.

Consider the relation $\Com(y)$ defined by
\[
\Com(y) \leftrightarrow y \in K(t) \land \exists \,x \in K(t): y^2 = x^3 -4.
\]
Since $y^2 = x^3 -3$ is a curve of genus 1, it does not admit a
rational parameterization, and so if an element $y$ satisfies
$\Com(y)$, then $y$ lies in $K$. Also, Denef~(\cite{Den78}) showed that
for every rational number $z$, there exists a rational number $y>z$
satisfying $\Com(y)$. We are now ready to define the set $T$.
\begin{theorem}\label{defineT}
Define the set $T$ by 
\begin{align}\nonumber
  Z \in T \leftrightarrow \;\;&\exists\, X_1, \dots, X_5,y \in K(t):\; \\\nonumber
  &\left(\Com(y)\; \land\right.\\ &\left.(y-t)Z^2 +1 = X_1^2 +X_2^2 +\dots
    + X_5^2\right).\label{order}
\end{align}
Then $T$ has the properties as in Remark~\ref{rem1}.
\end{theorem}
\begin{proof}
We follow the proof in \cite{Den78}: We will first show that every
element $Z \in T$ has positive order at $t^{-1}$.

Suppose there exist $X_1, \dots, X_5,y$ in $K(t)$ as in
Equation~(\ref{order}), and assume by contradiction that
$\ord_{t^{-1}}(Z) \leq 0$. Then $\deg Z \geq 0$, where $\deg Z$ denotes
the degree of the rational function $Z$. Since $y$ satisfies
$\Com(y)$, we have $y \in K$, which implies that $\deg(y-t) =1$, and so
$\deg((y-t)Z^2 +1)$ is positive and odd. But the degree of the
rational function $X_1^2 + \dots X_5^2$ is even, since in a formally
real field, a sum of squares is zero if and only if each term is zero,
and hence there is no cancellation of the coefficients of largest
degree.  Hence the left-hand-side of (\ref{order}) has odd degree,
while the right-hand-side has even degree, contradiction.

To show that the set $T$ satisfies the second property, let $Z \in
\Q(t)$, and assume $\ord_{t^{-1}}(Z) >0$. We want to show that $Z \in
T$. Since $\ord_{t^{-1}}(Z) >0$, we have $\ord_{t^{-1}}(tZ^2) >0$, and
so $tZ^2(r) \to 0$ as $|r| \to \infty$ ($r \in \R$). Hence we can find a
natural number $n$, such that for real numbers $r$ with $|r| >n$, we
have $|tZ^2 (r)| \leq 1/2.$ Pick a rational number $y$ with $y>n>0$ and
satisfying $\Com(y)$. Such a $y$ exists by the discussion before
Theorem~\ref{defineT}. Then
\[\left((y-t)Z^2 +1\right)(r) = yZ^2(r) -  tZ^2(r)+1 \geq yZ^2(r)-1/2+1
> 0\] for all $r \in \R$. By Pourchet's theorem, every positive
definite rational function over $\Q$ can be written as a sum of five
squares in $\Q(t)$. Hence there exist $X_1,\dots,X_5 \in K(t)$ as desired.
\end{proof}

\section{Function fields over the complex numbers in two
  variables}\label{KimandRoush}
Unfortunately, the diophantine definition of the set $T$ which defined
the elements of positive order at $t^{-1}$ and which was crucial for the
proof of Theorem~\ref{thm1} only works for formally real fields. 

For $\C(t_1,t_2)$ and finite extensions we will do something else that
avoids defining order.
\subsection{Hilbert's Tenth for the rational function field $\C(t_1,t_2)$}
In this section we will outline the proof of the following 
\begin{theorem}[\cite{KR92}]\label{TheoremKR}
Hilbert's Tenth Problem for $\C(t_1,t_2)$  with coefficients in
$\Z[t_1, t_2]$ is undecidable.
\end{theorem}
To prove undecidability of Hilbert's Tenth Problem for
$K:=\C(t_1,t_2)$ we will construct a diophantine model of the
structure $$\mathcal{S}:=\langle \Z \times \Z, + , \mid , \mathcal{Z},
\mathcal{W}\rangle $$ in $K$ (with coefficients in $\Z[t_1,t_2]$).  Here
$+$ denotes the usual com\-po\-nent-wise addition of pairs of integers,
$\mid$ re\-presents a relation which satisfies
\[(n,1)\mid(m,s)\Leftrightarrow m = ns,\]
and $\mathcal{Z}$ is a unary predicate which is interpreted as
\[
\mathcal{Z}(n,m)\Leftrightarrow m=0.
\]
The predicate $\W$ is interpreted as 
\[
\W((m,n),(r,s)) \Leftrightarrow m=s \land n=r.
\]

A {\em diophantine model} of $\mathcal{S}$ over $K$ is a diophantine
subset $S \subseteq K^n$ equipped with a bijection $\phi: \Z \times \Z \to S$ such
that under $\phi$, the graphs of addition, $\mid$, $\mathcal{Z}$, and
$\mathcal{W}$ in $\Z \times \Z$ correspond to diophantine subsets of $S^3,
S^2, S,$ and $S^2$, respectively.

A {\em diophantine model} of $\mathcal{S}$ over $K$ with coefficients
in $\Z[t_1,t_2]$ is a model, where in addition $S$ and the graphs of
addition, $\mid$, $\mathcal{Z}$, and $\mathcal{W}$ are diophantine over
$K$ with coefficients in $\Z[t_1,t_2]$.

We will now show that constructing such a
model is sufficient to prove undecidability of Hilbert's Tenth Problem
for $K$. First we can show the following
\begin{prop}$($\cite{Eis04}$)$
The relation $\W$ can be defined entirely in terms of the
other relations.
\end{prop}
\begin{proof}
It is enough to verify that
\[
\W((a,b),(x,y)) \Leftrightarrow (1,1)\mid((x,y)+(a,b)) \land
        (-1,1)\mid((x,y)-(a,b)).\]
\end{proof}

As Pheidas and Zahidi~(\cite{PheZa00}) point out we can existentially
define the integers with addition and multiplication inside
$$\mathcal{S}= \langle \Z \times
\Z,+,\mid\,,\mathcal{Z},\W\rangle,$$ so $\mathcal{S}$ has an
undecidable positive existential theory:

\begin{prop}
The structure $\mathcal{S}$ has an undecidable positive existential theory.
\end{prop}
\begin{proof}
  We interpret the integer $n$ as the pair $(n,0)$. The set $\{(n,0):n
  \in \Z\}$ is existentially definable in $\mathcal{S}$ through the
  relation $\mathcal{Z}$. Addition of integers $n,m$ corresponds to the
  addition of the pairs $(n,0)$ and $(m,0)$. To define multiplication
  of the integers $m$ and $r$, note that $n=mr$ if and only if
  $(m,1)\mid(n,r)$, hence $n=mr$ if and only if
\[
\exists\, a,b \,\,:((m,0)+(0,1))\mid((n,0)+(a,b))\land
\W((a,b),(r,0)).
\]
Since the positive existential theory of the integers with addition
and multiplication is undecidable, $\mathcal{S}$ has an undecidable positive
existential theory as well.
\end{proof}
The above proposition shows that in order to prove
Theorem~\ref{TheoremKR} it is enough to construct a diophantine model
of $\mathcal{S}$ over $K$ with coefficients in $\Z[t_1,t_2]$.  In the
next section we will construct this model.

\subsection{Generating elliptic curves of rank one}\label{elliptic}

As before, let $K:=\C(t_1,t_2)$.  Our first task is to find a
diophantine set $A$ over $K$ which is isomorphic to $\Z \times \Z$ as a
set.  Following Kim and Roush (\cite{KR92}) we will obtain such a set
by using the $K$-rational points on two elliptic curves which have
rank one over $K$.  The same argument as in Theorem~\ref{pointorder}
shows that the following proposition holds:
\begin{prop}
Let $E$ be an elliptic curve over $\Q$ without complex multiplication
and with Weierstrass equation $ y^2 = x^3 + ax +b$, where $a,b \in \Q$
and $b \neq 0$. Consider the twists $\EE_1, \EE_2$ of $E$ defined by
%
\[
\EE _1 :(t_1^3 + at_1 +b)Y^2 = X^3 + aX + b
\]
and 
\[
\EE_2: (t_2^3 + at_2 +b)Y^2 = X^3 + aX + b.
\]
The point $(t_i,1) \in \EE_i(K)$ has infinite order for
$i=1,2$, and $(t_i,1)$ generates $\EE_i(K)$ modulo points of
order 2.
\end{prop}

To be able to define a suitable set $S$ which is isomorphic to\linebreak $\Z \times
\Z$ we need to work in an algebraic extension $F$ of $K$.
Let $F:=\C(t_1,t_2)(h_1,h_2)$, where $h_i$ is defined by $h_i^2=t_i^3
+ at_i+b$, for $i=1,2$.

To prove that the positive existential theory of $K$ in the language
$\langle\,+,\cdot\;;0,1,t_1,t_2 \rangle$ is undecidable, it is enough to prove that
the positive existential theory of $F$ in the language $\langle\, +,\cdot \; ;
0,1,t_1,t_2,h_1,h_2,\mathcal{P}\,\rangle$ is undecidable, where
$\mathcal{P}$ is a predicate for the elements of the subfield $K$
(\cite[Lemma 1.9]{PheZa00}). So from now on we will work with equations
over $F$.

Over $F$ both $\EE_1$ and $\EE_2$ are isomorphic to $E$. There is an
isomorphism between $\EE_1$ and $E$ that sends $(x,y) \in \EE_1$ to the
point $(x,h_1y)$ on $E$. Under this isomorphism the point $(t_1,1)$ on
$\EE_1$ corresponds to the point $P_1:=(t_1,h_1)$ on $E$.  Similarly
there is an isomorphism between $\EE_2$ and $E$ that sends the point
$(t_2,1)$ on $\EE_2$ to the point $P_2:=(t_2,h_2)$ on $E$.

So the element $(n,m) \in \Z \times \Z$ corresponds to the point $nP_1 +
mP_2 \in E(F)$.
As in Section~\ref{obtainS}, we can take care of the point at $\infty$
on the curve $E$.

The set of points $\Z P_1 \times \Z P_2 \subseteq E(F)$ is existentially definable
in our language, because we have a predicate for the elements of
$K$ : Since $\EE_1$ has $2$-torsion, we first give a
diophantine definition of $2 \cdot \, \Z P_1$ as in
Section~\ref{obtainS}:
\[
P \in 2 \cdot \,\Z P_1 \Leftrightarrow \exists \,x,y \in K \,\, (t_1^3 +
a t_1 + b)\,y^2 = x^3 + ax + b \land P = 2 \cdot (x, h_1y)
\]
Then $\Z P_1$ can be defined as
\[
P \in \Z P_1 \Leftrightarrow (P \in 2\cdot \, \Z  P_1) {\mbox{ or }}
(\exists \,Q \in 2 \cdot \, \Z P_1\mbox{ and } P = Q + P_1)
\]
Similarly we have a diophantine definition for $\Z P_2$. Hence the
cartesian product $\Z P_1 \times \Z P_2 \subseteq E(F)$ is
existentially definable, since addition on $E$ is existentially
definable.

\subsection{Existential definition of $+$ and $\mathcal{Z}$}\label{relations}
The unary relation $\mathcal{Z}$ is existentially definable, since
this is the same as showing that the set $\Z P_1$ is diophantine,
which was done above.  Addition of pairs of integers corresponds to
addition on the cartesian product of the elliptic curves $\EE_i$ (as
groups), hence it is existentially definable.  Since $\W$ can be
defined in terms of the other relations, it remains to define the
divisibility relation $\mid$\,.
\subsection{Existential definition of divisibility}
In the following $x(P)$ will denote the $x$-coordinate of a
point $P$ on $E$, and $y(P)$ will denote the $y$-coordinate of
$P$. The following theorem gives the existential definition of $\mid$\,:
\begin{theorem}\label{definediv}
\begin{align*} &\forall m\in \Z,n,r \in \Z-\{0\}:\\ (m,1)\mid(n,r) &\Leftrightarrow\\ 
  &(\,\exists \, z,w \in F^{*}\, \,\, x(nP_1+rP_2)\, z^2 +
    x(mP_1+P_2)\,w^2 =1)
\end{align*}
\end{theorem}
Clearly this definition is existential in $(m,1)$ and $(n,r)$. It is
enough to give an existential definition of $\mid$ for $n,r \in
\Z-\{0\}$, because we can handle the cases when $n$ or $r$ are zero
separately.
\begin{proof}
 For the first implication, assume that $(m,1)\mid(n,r)$, {\it i.e.}\
  $n=mr$.  Then both $ x(nP_1+rP_2)=x(r(mP_1+P_2))$ and $x(mP_1+P_2)$
  are elements of $\C(x(mP_1+P_2), y(mP_1+P_2))$, which has
  transcendence degree one over $\C$. This means that we can apply the
  Tsen-Lang Theorem (Theorem~\ref{Tsen-Lang} from the appendix) to the
  quadratic form
\[
x(nP_1+rP_2)z^2 + x(mP_1+P_2)w^2 -v^2
\]
to conclude that there exists a nontrivial zero $(z,w,v)$ over
$\C(x(mP_1+P_2), y(mP_1+P_2))$. From the theory of quadratic forms it
follows that there exists a nontrivial zero $(z,w,v)$ with $z\cdot w
\cdot v\neq 0$.

For the other direction, suppose that $n \neq mr$ and assume by
contradiction that there exist $z,w
\in F^*$ with \begin{equation}\label{eq:1}x(nP_1+rP_2)\, z^2 +
    x(mP_1+P_2)\,w^2 =1.\end{equation}

\noindent{\bf Claim:} There exists a discrete valuation $w_m:F^{*}
\twoheadrightarrow \Z$ such that $w_m(x(mP_1+P_2))=1$ and such that
$w_m(x(nP_1+rP_2))=~0$.

\noindent{\bf Proof of Claim:} Let $P_2' = mP_1 + P_2
=(t_2',h_2')$. Remember that $F= \C(t_1,t_2,h_1,h_2)$. Then
$$F=\C(t_1,t_2,h_1,h_2)=\C(t_1,h_1,t_2',h_2')=
\C(x(P_1),y(P_1),x(P_2'),y(P_2')),$$ since $t_2 = x(P_2'-mP_1)$ and
$h_2=y(P_2'-mP_1)$.

Now let $w_m:F^{*}\twoheadrightarrow \Z$ be a
discrete valuation which extends the discrete valuation $w$ of
$\C(t_1,h_1)(t_2')$ associated to $t_2'$. The valuation $w$ is the
discrete valuation that satisfies $w(\gamma)=0$ for all $\gamma \in
\C(t_1,h_1)$ and $w(t_2')=1$.

Let $s:= n-mr$.  By assumption $s\neq 0$. We have $nP_1
+rP_2=sP_1+rP_2'$. The residue field of $w_m$ is $\C(t_1, h_1)$. Let
$x_{s,r}$ denote the image of $x(sP_1+rP_2')$ in the residue field of
$w_m$. Then $x_{s,r}=x\left(s(t_1,h_1)+r(0,\pm \sqrt b)\right)$.  We can
show that this $x$-coordinate cannot be zero, which will imply that
$w_m(x(nP_1+rP_2))=w_m(x(sP_1+rP_2')) =0$: The point $(t_1,h_1) \in
E(\C(t_1,h_1))$ has infinite order, $t_1$ is transcendental over $\C$,
and all points of $E$ whose $x$-coordinate is zero are defined over
$\C$. Since $s\neq 0$, this implies that
$x\left(s(t_1,h_1)+r(0,\pm \sqrt b)\right) \neq 0$. This proves
the claim.

Since $z,w$ satisfy Equation~(\ref{eq:1}) it easily follows that
$x_{s,r}$ is a square in the residue field.  This will give us a
contradiction:

We have $x_{s,r}=x\left(s (t_1,h_1) + r (0,\pm \sqrt{b})\right)$.
Since the residue field of $w_m$ is $\C(t_1,h_1)$, which is the
function field of $E$, we can consider $x_{s,r}$ as a function $E \to
\Proj^1_{\C}$. Then $x_{s,r}$ corresponds to the function on $E$ which
can be obtained as the composition $P \mapsto sP+ r (0, \sqrt
b)\mapsto x(sP+r(0,\sqrt b))$.  The $x$-coordinate map is of degree
$2$ and has two distinct zeros, namely $(0, \sqrt b)$ and $(0, -\sqrt
b)$.  The map $E \to E$ which maps $P$ to $(sP + r (0,\sqrt b))$ is
unramified since it is the multiplication-by-$s$ map followed by a
translation.  Hence the composition of these two maps has $2s^2$
simple zeros. In particular, it is not a square in $\C(t_1,h_1)$. This
completes the proof of the theorem.
\end{proof}

\section{Generalization to finite extensions of $\C(t_1,t_2)$}\label{general}
In \cite{Eis04} we proved the following theorem.
\begin{theorem}
  Let $L$ be the function field of a surface over the complex numbers.
  There exist $z_1, z_2 \in L$ that generate an extension of
  transcendence degree 2 of $\C$ and such that Hilbert's Tenth Problem
  for $L$ with coefficients in $\Z[z_1,z_2]$ is undecidable.
\end{theorem}
\begin{rem}
  This theorem also holds for transcendence degree $\geq 2$, {\it i.e.\
  }for finite extensions of $\C(t_1, \dots,t_n)$, with $n \geq 2$, and
  the proof of the more general theorem can also be found in
  \cite{Eis04}. To make our exposition as short as possible and to
  avoid extra notation, we will only discuss the transcendence degree
  2 case here.
\end{rem}

Our proof will proceed as the proof for $\C(t_1,t_2)$, {\it i.e.\ }we
will construct a diophantine model of $\langle \Z \times \Z, +
, \mid , \mathcal{Z}, \mathcal{W}\rangle $ in $L$.
We will now give an outline of the steps that are needed in the proof.

\subsection{Finding suitable elliptic curves of rank one}
In the above undecidability proof for $\C(t_1,t_2)$ (\cite{KR92}) we
obtained a set that was in bijection with $\Z \times \Z$ by using the
$\C(t_1,t_2)$-rational points on two elliptic curves which have rank
one over $\C(t_1,t_2)$. However, the two elliptic curves that we used
in Section~\ref{elliptic} could have a rank higher than one over
$L$, so we need to construct two new elliptic curves which
will have rank one over $L$.

To do this we use a theorem by Moret-Bailly (\cite[Theorem
1.8]{MB}):

\begin{theorem}\label{MB}
  Let $k$ be a field of characteristic zero. Let $C$ be a smooth
  projective geometrically connected curve over $k$ with function
  field $F$. Let $Q$ be a finite nonempty set of closed points of $C$.
  Let $E$ be an elliptic curve over $k$ with Weierstrass equation $y^2
  = x^3 + ax +b$ where $a,b \in k$ and $b \neq 0$.  Let $f \in F$ be
  admissible for $E,Q$.  For $\lambda \in k^{*}$ consider the twist $\EE_{\lambda
    f}$ of $E$ which is defined to be the smooth projective model of
  $$((\lambda f)^3 + a(\lambda f) +b)\,y^2 = x^3 + ax + b.$$ Then the natural
  homomorphism $\EE(k(T))\hookrightarrow \EE_{\lambda f}(F)$ induced by the
  inclusion $k(T) \hookrightarrow F$ that sends $T$ to $\lambda f$ is an
  isomorphism for infinitely many $\lambda \in \Z$.
\end{theorem}
We will not define here what it means to be admissible, but we will
only state that given $C,E,Q$ as above, admissible functions exists,
and if $f$ is admissible for $C,E,Q$, then for all but finitely many
$\lambda \in k^{*}$, $\lambda f$ is still admissible.

Now we can state our theorem that allows us to obtain two suitable
elliptic curves of rank one over $L$:

\begin{theorem}\label{rank2curves}
  Let $L$ be a finite extension of $\C(t_1,t_2)$.  Let $E/ \C$ be an
  elliptic curve with Weierstrass equation $$y^2 = x^3 + ax + b,$$
  where $a,b \in \C$, $b\neq 0$. Assume that $E$ does not have complex
  multiplication.  There exist $z_1,z_2 \in L$ such that $\C(z_1,z_2)$
  has transcendence degree $2$ over $\C$ and such that the two
  elliptic curves $\EE_1,\EE_2$ given by the affine equations $\EE_1:
  (z_1^3 + az_1+b)\,y^2 = x^3 + ax + b$ and $\EE_2: (z_2^3 +
  az_2+b)\,y^2 = x^3 + ax + b$ have rank one over $L$ with generators
  $(z_1,1)$ and $(z_2,1)$, respectively (modulo 2-torsion).
\end{theorem}
\begin{proof}
  This is proved in \cite{Eis04}. In the proof we apply
  Theorem~\ref{MB} with $k$ chosen to be the algebraic closure of
  $\C(t_2)$ inside $L$.
\end{proof}
\subsection{Existential definition of divisibility}
To existentially define the relation $\mid$, we need an elliptic curve
$E$ as in Theorem~\ref{rank2curves} with the additional property that
the point $(0,\sqrt{b})$ has infinite order. So from now on we fix
$E$ to be the smooth projective model of $y^2=x^3+x+1$. This curve
does not have complex multiplication, and the point $(0,1)$ has
infinite order (see the curve 496A1 in \cite{Cremona}). We fix
$z_1,z_2$ as in Theorem~\ref{rank2curves}.

As before, let $F:=\C(z_1,z_2)(h_1,h_2)$, where $h_i$ is defined by
$h_i^2=z_i^3 + az_i+b$, for $i=1,2$. Let $M:= L(h_1, h_2)$. Over $M$,
the elliptic curves $\EE_1$ and $\EE_2$ are isomorphic to $E$. Let
$P_1:=(z_1,h_1)$, $P_2:=(z_2,h_2)$ be the two points on $E$ as before.
From now on we will work with equations over $M$.  To give a
diophantine definition we would like to prove an analogue of
Theorem~\ref{definediv}. To make this theorem work we have to
introduce extra equations.

Let $\alpha := [M:F]$. We have the following theorem:
\begin{theorem}\label{almostdiv}\cite{Eis04}
  There exists a finite set $U \subseteq \Z$ such that for all $m \in
  \Z-U$ we have: for all $n,r \in \Z-\{0\}$
\begin{align*} (m,1)\mid(n,r) \Leftrightarrow\\ 
  \left( \,\exists \, y_0,z_0 \in M^{*}\, \,\, \right.&x(nP_1+rP_2)\, y_0^2 +
    x(mP_1+P_2)\,z_0^2 =1 \\ \land \,\, \exists\, y_1,z_1 \in
  M^{*} \,\,\,
  & x(2nP_1+2rP_2)\, y_1^2 + x(mP_1+P_2)\,z_1^2 =1 \\
  \cdots\\  \land \,\, \exists\,
  y_{\alpha},z_{\alpha} \in M^{*} \,\,\,
  & \left. x(2^{\alpha}nP_1+2^{\alpha}rP_2)\, y_{\alpha}^2 +
  x(mP_1+P_2)\,z_{\alpha}^2 =1 \, \right) \,.
\end{align*}
\end{theorem}
\begin{proof}[Outline of Proof]
  By the same argument as in the proof of Theorem~\ref{definediv}, if
  $n=mr$, then the $\alpha +1$ equations can all be satisfied.

  For the other direction, the exceptional set $U$ is necessary here
  because as in the proof for $\C(t_1,t_2)$, for each $m$, we
  construct a discrete valuation $w_m:M^{*} \twoheadrightarrow \Z$.
  This valuation $w_m$ extends a certain other discrete valuation
  $v_m: F^{*} \twoheadrightarrow \Z$. We have to exclude all integers
  $m$, for which $w_m | v_m$ is ramified, and we define $U$ to be this
  set of integers. Then $U$ is finite by Theorem~\ref{ramified} from the
  appendix.

Assume that $n \neq mr$, and let $s:=n-mr$. Assume by contradiction that
we can satisfy all $\alpha +1$ equations.  We can show that for all $m \in
\Z-U$ there exists a discrete valuation $w_m:M^{*} \twoheadrightarrow
\Z$ such that $w_m(x(mP_1+P_2))=1$ and such that
$w_m(x(knP_1+krP_2))=~0$ for $k=1,2,4,\dots,2^{\alpha}$.  Let
$P_2'=mP_1+P_2$, and denote by $x_{s,r}$ the image of
$x(sP_1+rP_2')=x(nP_1+rP_2)$ in the residue field $\ell$ of $w_m$. The
proof of Theorem~\ref{almostdiv} proceeds by first showing that the
elements $x_{s,r}\;, \;\dots,\; x_{2^{\alpha}s\;,\;2^{\alpha}r}$ are not
squares in $\C(z_1,h_1)$, and then by proving that the images of
$x_{s,r}, \dots, x_{2^{\alpha}s,2^{\alpha}r}$ in
\[V:=[(\ell^*)^2 \cap\, \C(z_1,h_1)^*]/(\C(z_1,h_1)^*)^2 \text{ are
  distinct.}\] But using Kummer theory one can show that since
$[\ell:\C(z_1,h_1)] \leq \alpha$, the size of $V$ is bounded by $\alpha$ as well.
This gives us the desired contradiction.\end{proof}

Once we have Theorem~\ref{almostdiv}, it is easy to define the
relation $\mid$ for
all $m \in \Z$ as follows:

  Let $m_0$ be a fixed element in $\Z-U$, and let $d$ be a positive
  integer such that $U \subseteq (m_0-d,m_0+d)$.
  Since \[n=mr \Leftrightarrow dn+m_0r=dmr+m_0r=(dm+m_0)r,\] we have
\[(m,1)\mid(n,r)
\Leftrightarrow (dm+m_0,1)\mid(dn+m_0r,r),\] and we can just work with
that formula instead. So
\begin{align*}(m,1)\mid&(n,r) \Leftrightarrow\\ &\exists a,b\: (dm+m_0,1)\mid((dn,r)+m_0(a,b)) \land
\mathcal{W}((a,b),(0,r)).\end{align*} It is an easy exercise to show that
$\mathcal{W}$ is existentially definable using
Theorem~\ref{almostdiv}.  Since $m_0$ is a fixed integer, this
together with Theorem~\ref{almostdiv} implies that the last expression
is existentially definable in $(m,1)$ and $(n,r)$.

\section{Appendix}\label{appendix}
In this section we state two theorems that we needed in
Sections~\ref{KimandRoush} and \ref{general}.
\begin{theorem}\label{ramified}
  Let $L$ and $K$ be function fields of one variable with constant
  fields $C_L$ and $C_K$, respectively, such that $L$ is an extension
  of $K$. If $L$ is separably algebraic over $K$, then there are at
  most a finite number of places of $L$ which are ramified over $K$.
\end{theorem}
\begin{proof}
This theorem is proved on p.~111 of \cite{Deu73} when $C_L \cap K
=C_K$, and the general theorem also follows.
\end{proof}
\begin{theorem}{Tsen-Lang Theorem.}\label{Tsen-Lang}
Let K be a function field of transcendence degree $j$ over an
algebraically closed field $k$. Let $f_1,\cdots,f_r$ be forms in $n$
variables over $K$, of degrees $d_1, \cdots, d_r$. If
\[
n > \sum_{i=1}^{r}d_i^j
\]
then the system $f_1=\cdots=f_r=0$ has a non-trivial zero in $K^n$.
\end{theorem}
\begin{proof}
  This is proved in Proposition 1.2 and Theorem 1.4 in Chapter 5 of
  \cite{Pfi95}.
\end{proof}

\end{document}